\documentclass[12pt]{article}
\usepackage{graphicx}

\usepackage{cite,amsfonts,amsmath}%,subeqn}
\usepackage{epic}

\def\be{\begin{equation}}
\def\eqn#1{\be\label{#1}}
\def\ee{\end{equation}}

\def\bea{\begin{eqnarray}}

\def\eea{\end{eqnarray}}

\newcommand{\eqna}[1]{\begin{subequations} \label{#1}
\begin{eqnarray}}
\def\eena{\end{eqnarray}
\end{subequations}}

 \def\ha{{\textstyle{\frac{1}{2}}}}

\def\ve{\varepsilon}

\def\bbc{C\kern-6.5pt I}
\def\ca{{\cal A}}
\def\nn{\nonumber}
\def\id{{\bf 1}}

\def\nt{\noindent}
\def\bu{$\bullet$~~~}
\def\lra{\longleftrightarrow}

\normalbaselines \oddsidemargin 0.5cm \evensidemargin 0.5cm
\begin{document}

\begin{center}
 \textsf{\LARGE Exotic Bialgebras from 9$\times$9 Unitary Braid Matrices}\footnote{Invited talk by VKD
 at XIII International Conference on Symmetry Methods in Physics, (Dubna, 2009), to appear in  Physics of Atomic
Nuclei.}

\vspace{10mm}
{\bf \large B. Abdesselam$^{a,}$\footnote{Email:  boucif@yahoo.fr}, A.
Chakrabarti$^{b,}$\footnote{Email:
chakra@cpht.polytechnique.fr},\\[2mm]
V.K. Dobrev$^{c,d,}$\footnote{Email: dobrev@inrne.bas.bg} and S.G.
Mihov$^{c,}$\footnote{Email: smikhov@inrne.bas.bg}}

\vspace{4mm}
  \emph{$^a$ Laboratoire de Physique Quantique de la
Mati\`ere et de Mod\'elisations Math\'ematiques, Centre
Universitaire de Mascara, 29000-Mascara, Alg\'erie\\
and \\
Laboratoire de Physique Th\'eorique, Universit\'e d'Oran
Es-S\'enia, 31100-Oran, Alg\'erie}
  \\
  \vspace{2mm}
  \emph{$^b$ Centre de Physique Th{\'e}orique, Ecole Polytechnique, 91128 Palaiseau Cedex, France.}
  \\
  \vspace{2mm}
  \emph{$^c$ Institute of Nuclear Research and Nuclear Energy
Bulgarian Academy of Sciences 72 Tsarigradsko Chaussee, 1784
Sofia, Bulgaria}\\
 \vspace{2mm}
 \emph{$^d$ Abdus Salam International Center for Theoretical Physics
Strada Costiera 11, 34100 Trieste, Italy}
\end{center}

\vskip 1cm

\begin{abstract}
We construct the exotic bialgebras that arise from a $9\times9$ unitary braid matrix.
 We also construct the dual bialgebra of one of these exotic bialgebras.
\end{abstract}

\vskip 2cm

\section{Introduction}
\label{sect:intro}
\setcounter{equation}{0}

For several years \cite{ACDM1,ACDM2,ACDM3,ACDM4}
our initial collaboration studied the   algebraic structures coming from
4x4 $R$-matrices (solutions of the Yang--Baxter equation)
that are not deformations of classical ones (i.e., the
identity up to signs). In parallel, higher
dimensional ($N^2 \times N^2$ matrices for all $N>2$) exotic braid matrices   have been presented and studied in
\cite{AC0401,Chakra3}. More recently, our follow-up collaboration (with Boucif Abdesselam replacing our deceased
friend and co-author Daniel Arnaudon) constructed $N^2\times N^2$ unitary braid
matrices $\widehat{R}$ for $N>2$ generalizing the class known
for $N=2$ \cite{ACDM5,ACDM6}. Some of these results were applied to model-building \cite{AbC1,AbC2}.

In the present paper we start the systematic study of the bialgebras
that arise from these higher dimensional unitary braid matrices. We
start with the simplest possible case ~$N=3$~ in order to get the
necessary expertise. However, even this case is complicated enough.

The paper is organized as follows. In Section 2 we give a the
overall general setting. Section  3 is  devoted to the
 explicit construction of  exotic bialgebras starting from a $9\times9$   braid matrix
that we constructed in \cite{ACDM5}. In Section  4   we
 find the dual bialgebra  of the most promising exotic bialgebra
 constructed in Section 3.

\section{Preliminaries}
\setcounter{equation}{0}

Our starting point is the following $9\times9$   braid matrix from
\cite{ACDM5}:
\begin{equation}\label{braidnine}
\widehat{R}(\theta)=\left|\begin{array}{ccccccccc}
   a_+ & 0 & 0 & 0 & 0 & 0 & 0 & 0 & a_-\\
   0 & b_+ & 0 & 0 & 0 & 0 & 0 & b_- & 0\\
   0 & 0 & a_+ & 0 & 0 & 0 & a_- & 0 & 0\\
   0 & 0 & 0 & c_+ & 0 & c_- & 0 & 0 & 0\\
   0 & 0 & 0 & 0 & 1 & 0 & 0 & 0 & 0\\
   0 & 0 & 0 & c_- & 0 & c_+ & 0 & 0 & 0\\
   0 & 0 & a_- & 0 & 0 & 0 & a+ & 0 & 0\\
   0 & b_- & 0 & 0 & 0 & 0 & 0 & b_+ & 0\\
   a_- & 0 & 0 & 0 & 0 & 0 & 0 & 0 & a_+\\
   \end{array}\right|
\end{equation}
where
\begin{equation}\label{param}
a_\pm = \ha (e^{m_{11}^+ \theta} \pm e^{m_{11}^-
\theta}), \ \ b_\pm = \ha (e^{m_{21}^+ \theta} \pm e^{m_{21}^-
\theta}), \ \ c_\pm = \ha (e^{m_{22}^+ \theta} \pm e^{m_{22}^-
\theta}), \end{equation}
and ~$m^\pm_{ij}$~ are parameters. The above braid matrix satisfies the baxterized braid equation:
\begin{equation}
\widehat{R}_{12}(\theta)\widehat{R}_{23}(\theta+\theta')\widehat{R}_{12}(\theta')=
\widehat{R}_{23}(\theta')\widehat{R}_{12}(\theta+\theta')\widehat{R}_{23}(\theta).
\end{equation}

For the RTT relations of Faddeev-Reshetikhin-Takhtajan \cite{FRT} we
need the corresponding R-matrix, ~$R=P\hat{R}$, ($P$ is the
permutation matrix):
\begin{equation}\label{ninexnine}
R(\theta)=
\left(
\begin{array}{ccccccccc}
a_+ & 0 & 0 & 0 & 0 & 0 & 0 & 0 & a_- \cr
0 & 0 & 0 & c_+ & 0 & c_- & 0 & 0 & 0 \cr
0 & 0 & a_- & 0 & 0 & 0 & a_+ & 0 & 0 \cr
\ 0 & b_+ & 0 & 0 & 0 & 0 & 0 & b_- & 0 \cr
0 & 0 & 0 & 0 & 1 & 0 & 0 & 0 & 0 \cr
0 & b_- & 0 & 0 & 0 & 0 & 0 & b_+ & 0 \cr
0 & 0 & a_+ & 0 & 0 & 0 & a_- & 0 & 0 \cr
0 & 0 & 0 & c_- & 0 & c_+ & 0 & 0 & 0 \cr
a_- & 0 & 0 & 0 & 0 & 0 & 0 & 0 & a_+
\end{array}
\right)
\end{equation}
 which satisfies the baxterized Yang-Baxter equation:
\begin{equation}
{R}_{12}(\theta){R}_{13}(\theta+\theta'){R}_{23}(\theta')=
{R}_{23}(\theta'){R}_{13}(\theta+\theta'){R}_{12}(\theta)
\end{equation}

In fact, we need the solutions of the constant YBE which are as follows:
\begin{equation}\label{konst} a_+ = b_+ = c_+ = 1/2 , \ \ \ a_+ = \pm a_-, \ \ \ b_+ =
\pm b_-, \ \ \ c_+ = \pm c_- \end{equation}

In view of  \eqref{param} we see that for $a_+ = a_-$ the proper
limit is obtained, for example, by taking the following limits:
first ~$m_{11}^- = - \infty$, and then ~$\theta =0$, while for
$a_+ = - a_-$ the limit may be obtained for $m_{11}^+ = - \infty$ first,
and then $\theta = 0$. Similarly are obtained the limits for $b_\pm$ and
$c_\pm$.

So we have eight $R$ matrices satisfying the constant YBE~:
\eqn{varr}
 (+,+,+), \ (-,+,+), \ (+,-,+), \ (+,+,-), \ (+,-,-), \ (-,+,-),
\ (-,-,+), \ (-,-,-)\ee
where the $\pm$ signs denote respectively the signs of $a_+ = \pm
a_-$, $b_+ = \pm b_-$ and $c_+ = \pm c_-$.

For the elements of the $3\times3$ ~$T$ matrix we introduce the
notation:
\begin{equation}\label{mm}
 T = \begin {pmatrix} k & p & l \\ q & r & s
\\ m & t & n \\ \end{pmatrix} \end{equation}

\section{Solutions of the RTT equations and exotic bialgebras}
\setcounter{equation}{0}

 We consider  matrix bialgebras which are
unital associative algebras generated by the nine elements from
\eqref{mm} ~$k,l,m,n,p,q,r,s,t$. The co-product and co-unit
relations are the classical ones: \eqna{coal}
&& \delta \left( T \right)=   T   ~\otimes~   T   \\
&& \ve \left( T \right) = \id_3
\eena
\nt  We expect the bialgebras under consideration not to be Hopf
algebras which, as in the $S03$ case \cite{ACDM2}, would be easier to check after
we find the dual bialgebras.

In the next subsections we obtain the desired bialgebras  by applying the RTT relations of \cite{FRT}:
\eqn{rtt} R\ T_1\ T_2 \ \ =\ \ T_2\ T_1\ R \ \ , \end{equation}
where \ $T_1 \ =\ T\, \otimes\, \id_2$\ , \ $T_2 \ =\ \id_2 \,
\otimes\, T$, ~for   $\ R\ =\ R(\theta)\ $, \eqref{ninexnine},
and the parameters are the constants in \eqref{konst} following  the eight cases of \eqref{varr}.

\subsection{Algebraic relations}

I) Relations, which do not depend on the parameters $a_\pm, \
b_\pm, c_\pm$. We have the set of relations

\begin{eqnarray} (N) = \{ k^2 = n^2, \ \ kn = nk, \ \ l^2 = m^2, \ \ lm =
ml \nonumber \\ km = nl, \ \ kl = nm, \ \ lk = mn, \ \ mk = ln \nonumber \\
r(k-n) = (k-n)r = 0, \ \ r(l-m) = (l-m)r = 0 \} \end{eqnarray}

The last two relations suggest to introduce the generators:
\begin{eqnarray} k = \tilde{k} + \tilde{n}, \ \ n = \tilde{k} -
\tilde{n}; \ \  l = \tilde{l} + \tilde{m}, \ \ m = \tilde{l} -
\tilde{m}, \nonumber \\ p = \tilde{p} + \tilde{t}, \ \ t =
\tilde{p} - \tilde{t}; \ \  q = \tilde{q} + \tilde{s}, \ \ s =
\tilde{q} - \tilde{s}, \end{eqnarray}

In terms of these generators we have:
\begin{eqnarray} (N) = \{\tilde{k}\tilde{n} = \tilde{n}\tilde{k} =
0; \ \
\tilde{l}\tilde{m} = \tilde{m}\tilde{l} = 0; \nonumber \\
\tilde{k}\tilde{m} = \tilde{n}\tilde{l} = 0; \ \
\tilde{l}\tilde{n} = \tilde{m}\tilde{k} = 0; \nonumber \\
r\tilde{m} = r\tilde{n} = 0, \ \ \tilde{m}r = \tilde{n}r = 0 \}
\end{eqnarray}

II) Relations that do not depend on the relative signs of
$(a_-,b_-)$, $(a_-,c_-)$ and $(b_-, c_-)$. In that case we have

IIa) $a_+ = \pm a_-$

If $a_+ = a_-$ we have the set of relations
\begin{equation} A_+ = \{p^2 = t^2, \ \ pt = tp;  \ \ q^2 = s^2, \ \ qs
= sq \} \end{equation}
or in terms of the alternative generators we have:
\begin{equation} A_+ = \{\tilde{p}\tilde{t} = \tilde{t}\tilde{p} =
0; \ \ \tilde{q}\tilde{s} = \tilde{s}\tilde{q} = 0 \}
\end{equation}

If $a+ = - a_-$ the set of relations is
\begin{equation} A_- = \{p^2 = - t^2, \ \ pt = - tp; \ \ q^2 = - s^2, \ \ qs
= - sq \} \end{equation}
or alternatively
\begin{equation} A_- = \{\tilde{p}^2 = \tilde{t}^2= 0; \ \ \tilde{q}^2 = \tilde{s}^2 = 0
\} \end{equation}

IIb) $b_+ = \pm b_-$

If $b_+ = b_-$ we have the set of relations:
\begin{equation} B_+ = \{ rp = rt; \ \ rq = rs \} \end{equation}
Alternatively
\begin{equation} B_+ = \{ r\tilde{t} = 0; \ \ r\tilde{s} = 0  \} \end{equation}

If $b+ = - b_-$
\begin{equation} B_- = \{ rp = - rt; \ \ rq = - rs \} \end{equation}
Alternatively
\begin{equation} B_- = \{ r\tilde{p} = 0; \ \ r\tilde{q} = 0  \} \end{equation}

IIc) $c_+ = c_\pm$

If $c_+ = c_-$ we have
\begin{equation} C_+ = \{ pr = tr; \ \ qr = sr \} \end{equation}
Alternatively
\begin{equation} C_+ = \{ \tilde{t}r = 0; \ \ \tilde{s}r = 0  \} \end{equation}

If $c_+ = - c_-$ we have
\begin{equation} C_- = \{ pr = - tr; \ \ qr = - sr \} \end{equation}
Alternatively
\begin{equation} C_- = \{ \tilde{p}r = 0; \ \ \tilde{q}r = 0  \} \end{equation}

III) Relations depending on the relative signs of $(a_-, b_-)$,
$(a_-,b_-0)$ and $(b_-,c_-)$

IIIa) $a_- = \pm b_-$

If $a_- = b_-$ we have the set of relations
\begin{eqnarray} (AB)_+ = \{ pk = tn, \ tk = pn; \ pl = tm,  \ tl =
pm; \nonumber \\ qk = sn, \ qn = sk; \ ql = sm, \ sl = qm \}
\end{eqnarray}
Alternatively
\begin{equation} (AB)_+ = \{ \tilde{p}\tilde{m} =
\tilde{p}\tilde{n} = \tilde{t}\tilde{k} = \tilde{t}\tilde{l} = 0;
\ \ \tilde{q}\tilde{m} = \tilde{q}\tilde{n} = \tilde{s}\tilde{k} =
\tilde{s}\tilde{l} = 0 \} \end{equation}

If $a_- = - b_-$ we have
\begin{eqnarray} (AB)_- = \{ pk = - tn, \ tk = - pn; \ pl = - tm, \ tl
= - pm; \nonumber \\ qk = - sn, \ qn = - sk; \ ql = - sm, \ sl = -
qm \}
\end{eqnarray}
Alternatively
\begin{equation} (AB)_- = \{ \tilde{p}\tilde{k} =
\tilde{p}\tilde{l} = \tilde{t}\tilde{m} = \tilde{t}\tilde{n} = 0;
\ \ \tilde{q}\tilde{k} = \tilde{q}\tilde{l} = \tilde{s}\tilde{m} =
\tilde{s}\tilde{n} = 0 \} \end{equation}

IIIb) $a_- = \pm c_-$

If $a_- = c_-$ we have
\begin{eqnarray} (AC)_+ = \{ kp = nt, \ kt = np; \ lp = mt,  \ lt =
mp; \nonumber \\ kq = ns, \ nq = ks; \ lq = ms, \ ls = mq \}
\end{eqnarray}
Alternatively
\begin{equation} (AC)_+ = \{ \tilde{k}\tilde{t} =
\tilde{l}\tilde{t} = \tilde{m}\tilde{p} = \tilde{n}\tilde{p} = 0;
\ \ \tilde{k}\tilde{s} = \tilde{l}\tilde{s} = \tilde{m}\tilde{q} =
\tilde{n}\tilde{q} = 0\} \end{equation}

If $a_- = - c_-$ we have
\begin{eqnarray} (AC)_- = \{ kp = - nt, \ kt = - np; \ lp = - mt, \ lt
= - mp; \nonumber \\ kq = - ns, \ nq = - ks; \ lq = - ms, \ ls = -
mq \}
\end{eqnarray}
Alternatively
\begin{equation} (AC)_- = \{ \tilde{k}\tilde{p} =
\tilde{l}\tilde{p} = \tilde{m}\tilde{t} = \tilde{n}\tilde{t} = 0;
\ \ \tilde{k}\tilde{q} = \tilde{l}\tilde{q} = \tilde{m}\tilde{s} =
\tilde{n}\tilde{s} = 0\} \end{equation}

IIIc) $b_- = \pm c_-$

If $b- = c_-$ we have
\begin{equation} (BC)_+ = \{ pq = ts, \ tq = ps; \ \ qp = st, \ qt
= sp \} \end{equation}
Alternatively
\begin{equation} (BC)_+ = \{ \tilde{p}\tilde{s} =
\tilde{t}\tilde{q} = 0;  \ \ \tilde{s}\tilde{p} =
\tilde{q}\tilde{t} = 0\} \end{equation}

If $b_- = - c_-$ we have
\begin{equation} (BC)_- = \{ pq = - ts, \ tq = - ps; \ \ qp = - st, \ qt
= - sp\} \end{equation}
Alternatively
\begin{equation} (BC)_- = \{ \tilde{p}\tilde{q} =
\tilde{t}\tilde{s} = 0;  \ \ \tilde{q}\tilde{p} =
\tilde{s}\tilde{t} = 0\} \end{equation}

\subsection{Classification of bialgebras}

Thus we have the following solutions

\bigskip\nt\bu for $a_+ = a_- = b_- = c_-$ we have the set of relations
\begin{equation} (+,+,+) = \{ N \cup A_+ \cup B_+ \cup C_+ \cup
(AB)_+ \cup (AC)_+ \cup (BC)_+ \} \end{equation}

Explicitly we have:
\begin{eqnarray}\label{ppp}
\tilde{k}\tilde{m} = \tilde{m}\tilde{k} = 0;\ \
\tilde{k}\tilde{n} = \tilde{n}\tilde{k} = 0; \ \
\tilde{k}\tilde{t} = \tilde{t}\tilde{k} = 0; \ \
\tilde{k}\tilde{s} =\tilde{s}\tilde{k} = 0 ; \nn\\
\tilde{l}\tilde{m} = \tilde{m}\tilde{l} = 0;\ \
\tilde{l}\tilde{n}  = \tilde{n}\tilde{l} = 0; \ \
 \tilde{l}\tilde{t} = \tilde{t}\tilde{l} = 0;\ \
\tilde{l}\tilde{s} = \tilde{s}\tilde{l} = 0; \nn\\
\tilde{p}\tilde{m} =  \tilde{m}\tilde{p} = 0;\ \
\tilde{p}\tilde{n} = \tilde{n}\tilde{p} = 0;\ \
\tilde{p}\tilde{t} = \tilde{t}\tilde{p} = 0; \ \
\tilde{p}\tilde{s} = \tilde{s}\tilde{p} = 0; \nn\\
\tilde{q}\tilde{m} = \tilde{m}\tilde{q} =0;\ \
\tilde{q}\tilde{n} = \tilde{n}\tilde{q} = 0;\ \
\tilde{q}\tilde{t} = \tilde{t}\tilde{q} = 0;  \ \
\tilde{q}\tilde{s} = \tilde{s}\tilde{q} = 0; \nn\\
r\tilde{m} = \tilde{m}r = 0; \ \
r\tilde{n} =   \tilde{n}r = 0;\ \
r\tilde{t} = \tilde{t}r = 0; \ \
r\tilde{s} =  \tilde{s}r = 0.    \end{eqnarray}

From \eqref{ppp} we see that the algebra ~$\ca_{+++}$~ is a direct sum of two subalgebras:
~$\ca^1_{+++}$~ with generators ~$\tilde{k},\tilde{l},\tilde{p},\tilde{q},r$, and
~$\ca^2_{+++}$~ with generators ~$\tilde{m},\tilde{n},\tilde{s},\tilde{t}$. Both subalgebras are free, with
no relations, and thus, no PBW bases.

\bigskip\nt\bu for $a_+ = - a_- = b_- = c_-$ we have the set of relations
\begin{equation} (-,+,+) = \{ N \cup A_- \cup B_+ \cup C_+ \cup
(AB)_- \cup (AC)_- \cup (BC)_+ \} \end{equation}

Explicitly we have:
\begin{eqnarray}\label{mpp}
\tilde{k}\tilde{m} =\tilde{m}\tilde{k} = 0;\ \
\tilde{k}\tilde{n} = \tilde{n}\tilde{k} =0; \ \
\tilde{k}\tilde{p} = \tilde{p}\tilde{k} =0;\ \
\tilde{k}\tilde{q} = \tilde{q}\tilde{k} = 0;\nn\\
\tilde{l}\tilde{m} = \tilde{m}\tilde{l} = 0;  \ \
\tilde{l}\tilde{n} = \tilde{n}\tilde{l} = 0; \ \
\tilde{l}\tilde{p} = \tilde{p}\tilde{l} = 0;\ \
\tilde{l}\tilde{q} = \tilde{q}\tilde{l} = 0;\nn\\
r\tilde{m} = \tilde{m}r = 0;\ \
r\tilde{n} =  \tilde{n}r = 0;  \ \
r\tilde{t}= \tilde{t}r = 0;\ \
r\tilde{s} =   \tilde{s}r = 0; \nn\\
 \tilde{t}\tilde{m} = \tilde{m}\tilde{t} = 0;\ \
\tilde{t}\tilde{n} = \tilde{n}\tilde{t} = 0; \ \
\tilde{t}\tilde{q} = \tilde{q}\tilde{t} = 0;\nn\\
\tilde{s}\tilde{m} =\tilde{m}\tilde{s} = 0;\ \
\tilde{s}\tilde{n} = \tilde{n}\tilde{s} = 0;\ \
\tilde{s}\tilde{p} = \tilde{p}\tilde{s} = 0;\nn\\
\tilde{p}^2 = \tilde{q}^2=   \tilde{s}^2 = \tilde{t}^2 = 0.  \end{eqnarray}

The structure of this algebra, denoted ~$\ca_{-++}\,$, ~is more complicated. There are two quasi-free subalgebras:
~$\ca^1_{-++}$~ with generators ~$\tilde{k},\tilde{l},\tilde{t},\tilde{s}$, and
~$\ca^2_{-++}$~ with generators ~$\tilde{m},\tilde{n},\tilde{p},\tilde{q}$. They are quasi-free due to the last line of \eqref{mpp}.
They do not form a direct sum due to the existence of the following 12 two-letter building blocks of the basis ~$\ca_{-++}$~:
~$r\tilde{k}\,,r\tilde{l}\,,r\tilde{p}\,,r\tilde{q}\,,\tilde{p}\tilde{t}\,,\tilde{q}\tilde{s}$~ plus the reverse order.

\bigskip\nt\bu for $a_+ = a_- = - b_- = c_-$ we have
\begin{equation} (+,-,+) = \{ N \cup A_+ \cup B_- \cup C_+ \cup
(AB)_- \cup (AC)_+ \cup (BC)_- \} \end{equation}

Explicitly we have:
\begin{eqnarray}\label{pmp}
\tilde{k}\tilde{m} =\tilde{m}\tilde{k} =0;\ \
\tilde{k}\tilde{n} = \tilde{n}\tilde{k} =0; \ \
\tilde{k}\tilde{t} = \tilde{k}\tilde{s} = 0;\ \
\tilde{p}\tilde{k} = \tilde{q}\tilde{k} = 0;\nn\\
\tilde{l}\tilde{m} = \tilde{m}\tilde{l} = 0;  \ \
\tilde{l}\tilde{n} = \tilde{n}\tilde{l} = 0; \ \
\tilde{l}\tilde{t} = \tilde{l}\tilde{s} = 0; \ \
\tilde{p}\tilde{l} = \tilde{q}\tilde{l} = 0;\nn\\
r\tilde{m} = \tilde{m}r = 0;\ \
r\tilde{n} =   \tilde{n}r = 0; \ \
r\tilde{p} = r\tilde{q} = 0;\ \
\tilde{t}r = \tilde{s}r = 0; \nn\\
\tilde{t}\tilde{p} =\tilde{p}\tilde{t} =  0; \ \
\tilde{t}\tilde{s} = \tilde{s}\tilde{t} = 0; \ \
\tilde{t}\tilde{m} =\tilde{t}\tilde{n} = 0;\nn\\
\tilde{q}\tilde{p} = \tilde{p}\tilde{q} = 0;  \ \
\tilde{q}\tilde{s} = \tilde{s}\tilde{q} = 0; \ \
\tilde{m}\tilde{q} = \tilde{n}\tilde{q} = 0; \nn\\
 \tilde{s}\tilde{m} = \tilde{s}\tilde{n} = 0;\ \
 \tilde{m}\tilde{p} = \tilde{n}\tilde{p} = 0
.  \end{eqnarray}

The structure of this algebra, denoted ~$\ca_{+-+}\,$, ~is also complicated. There are four free subalgebras:
~$\ca^1_{+-+}$~ with generators ~$\tilde{k},\tilde{l},r$,
~$\ca^2_{+-+}$~ with generators ~$\tilde{m},\tilde{n}$,
~$\ca^3_{+-+}$~ with generators ~$\tilde{p},\tilde{s}$,
~$\ca^4_{+-+}$~ with generators ~$\tilde{q},\tilde{t}$. Only the first two are in direct sum,
otherwise all are related by the following 20 two-letter building blocks:
~$\tilde{k}\tilde{p}\,,\tilde{k}\tilde{q}\,,\tilde{s}\tilde{k},\,\tilde{t}\tilde{k}$,
~$\tilde{l}\tilde{p}\,,\tilde{l}\tilde{q}\,,\tilde{s}\tilde{l},\,\tilde{t}\tilde{l}$,
~$r\tilde{s}\,, r \tilde{t}\,, \tilde{p}r\,,  \tilde{q}r$,
~$\tilde{m}\tilde{s}\,, \tilde{m}\tilde{t}\,,
\tilde{p}\tilde{m}\,, \tilde{q}\tilde{m}$,
~$\tilde{n}\tilde{s}\,, \tilde{n}\tilde{t}\,,
\tilde{p}\tilde{n}\,, \tilde{q}\tilde{n}$.

There is no overall ordering. There is some partial order if we consider the subalgebra formed by
the generators of the latter three subalgebras: ~$\tilde{m},\tilde{n}$, ~$\tilde{p},\tilde{s}$,
~$\tilde{q},\tilde{t}$, namely, we have:
\eqn{ord} \tilde{p},\tilde{q} ~>~ \tilde{m},\tilde{n} ~>~ \tilde{s},\tilde{t} \ee
But for the natural subalgebra formed by generators ~$\tilde{k},\tilde{l},r$, ~$\tilde{p},\tilde{s}$,
~$\tilde{q},\tilde{t}$,  we have cyclic ordering:
\eqn{nord} \tilde{p},\tilde{q} ~>~ r ~>~ \tilde{s},\tilde{t} ~>~ \tilde{k},\tilde{l} ~>~ \tilde{p},\tilde{q} \ ,\ee
i.e., no ordering.
We have seen this phenomenon in the simpler exotic bialgebra ~$S03$, \cite{ACDM2}.

\bigskip\nt\bu for $a_+ = a_- = b_- = - c_-$ we have the set of relations

\begin{equation} (+,+,-) = \{ N \cup A_+ \cup B_+ \cup C_- \cup
(AB)_+ \cup (AC)_- \cup (BC)_- \} \end{equation}

Explicitly we have:
\begin{eqnarray}
\tilde{k}\tilde{m} = \tilde{m}\tilde{k} = 0; \ \
\tilde{k}\tilde{n} = \tilde{n}\tilde{k} =0; \ \
\tilde{k}\tilde{p} = \tilde{k}\tilde{q} = 0;\ \
\tilde{s}\tilde{k} = \tilde{t}\tilde{k} = 0;\nn\\
\tilde{l}\tilde{m} = \tilde{m}\tilde{l} = 0;  \ \
\tilde{l}\tilde{n} = \tilde{n}\tilde{l} = 0; \ \
\tilde{l}\tilde{p} = \tilde{l}\tilde{q} = 0;\ \
\tilde{s}\tilde{l} = \tilde{t}\tilde{l} = 0;\nn\\
r\tilde{m} = \tilde{m}r = 0;\ \
r\tilde{n} =    \tilde{n}r = 0; \ \
r\tilde{t} = r\tilde{s} = 0; \ \
 \tilde{p}r = \tilde{q}r = 0; \nn\\
\tilde{t}\tilde{p} =\tilde{p}\tilde{t} =  0; \ \
\tilde{t}\tilde{s} =\tilde{s}\tilde{t} = 0; \ \
\tilde{m}\tilde{t} =\tilde{n}\tilde{t} = 0;\nn\\
\tilde{q}\tilde{p} =\tilde{p}\tilde{q} = 0;  \ \
\tilde{q}\tilde{s} = \tilde{s}\tilde{q} = 0; \ \
\tilde{q}\tilde{m} =\tilde{q}\tilde{n} =   0;\nn\\
\tilde{p}\tilde{m} = \tilde{p}\tilde{n} =0; \ \
  \tilde{m}\tilde{s} = \tilde{n}\tilde{s} = 0. \end{eqnarray}

This algebra, denoted ~$\ca_{++-}\,$, ~is a conjugate of the previous one. It has the same four free algebras,
and the only difference is that the subalgebras are related by 20 two-letter building blocks which are
in reverse order w.r.t. the previous case:
~$\tilde{k}\tilde{s}\,,\tilde{k}\tilde{t}\,,\tilde{p}\tilde{k},\,\tilde{q}\tilde{k}$,
~$\tilde{l}\tilde{s}\,,\tilde{l}\tilde{t}\,,\tilde{p}\tilde{l},\,\tilde{q}\tilde{l}$,
~$r\tilde{p}\,, r \tilde{q}\,, \tilde{s}r\,,  \tilde{s}r$,
~$\tilde{m}\tilde{p}\,, \tilde{m}\tilde{q}\,,
\tilde{s}\tilde{m}\,, \tilde{t}\tilde{m}$,
~$\tilde{n}\tilde{p}\,, \tilde{n}\tilde{q}\,,
\tilde{s}\tilde{n}\,, \tilde{t}\tilde{n}$.

\bigskip\nt\bu for $a_+ = a_- = - b_- = - c_-$ we have the set of relations
\begin{equation} (+,-,-) = \{ N \cup A_+ \cup B_- \cup C_- \cup
(AB)_- \cup (AC)_- \cup (BC)_+ \} \end{equation}

Explicitly we have:
\begin{eqnarray}
\tilde{k}\tilde{m} =\tilde{m}\tilde{k} = 0;\ \
\tilde{k}\tilde{n} = \tilde{n}\tilde{k} =0; \ \
\tilde{k}\tilde{p} =\tilde{p}\tilde{k} = 0;\ \
\tilde{k}\tilde{q} = \tilde{q}\tilde{k} = 0;\nn\\
\tilde{l}\tilde{m} = \tilde{m}\tilde{l} = 0;  \ \
\tilde{l}\tilde{n} =  \tilde{n}\tilde{l} = 0; \ \
\tilde{l}\tilde{p} = \tilde{p}\tilde{l} = 0;\ \
\tilde{l}\tilde{q} = \tilde{q}\tilde{l} = 0;\nn\\
r\tilde{m} = \tilde{m}r = 0; \ \
r\tilde{n} =   \tilde{n}r = 0; \ \
r\tilde{p} = \tilde{p}r = 0;\ \
r\tilde{q} =   \tilde{q}r = 0; \nn\\
\tilde{s}\tilde{m} = \tilde{m}\tilde{s} = 0;\ \
\tilde{s}\tilde{n} = \tilde{n}\tilde{s} = 0; \ \
\tilde{s}\tilde{p} = \tilde{p}\tilde{s} = 0;\ \
\tilde{s}\tilde{q} = \tilde{q}\tilde{s} =  0; \nn\\
\tilde{t}\tilde{m} = \tilde{m}\tilde{t} = 0; \ \
\tilde{t}\tilde{n} = \tilde{n}\tilde{t} = 0; \ \
\tilde{t}\tilde{p} = \tilde{p}\tilde{t} =  0; \ \
\tilde{t}\tilde{q} = \tilde{q}\tilde{t} = 0.
\end{eqnarray}

This algebra, denoted ~$\ca_{+--}\,$, ~is a conjugate of the algebra ~$\ca_{+++}\,$, obtained
by the exchange of the pairs of generators ~$(\tilde{p},\tilde{q})$~ and ~$(\tilde{s},\tilde{t})$.

\bigskip\nt\bu for $a_+ = - a_- = b_- = - c_-$ we have the set of relations
\begin{equation} (-,+,-) = \{ N \cup A_- \cup B_+ \cup C_- \cup
(AB)_- \cup (AC)_+ \cup (BC)_- \} \end{equation}

Explicitly we have:
\begin{eqnarray}\label{mpm}
\tilde{k}\tilde{m} =\tilde{m}\tilde{k} = 0;\ \
\tilde{k}\tilde{n} = \tilde{n}\tilde{k} =0; \ \
\tilde{k}\tilde{t} =\tilde{k}\tilde{s} = 0; \ \
\tilde{p}\tilde{k} = \tilde{q}\tilde{k} = 0;\nn\\
\tilde{l}\tilde{m} = \tilde{m}\tilde{l} = 0;  \ \
\tilde{l}\tilde{n} =  \tilde{n}\tilde{l} = 0; \ \
\tilde{l}\tilde{t} = \tilde{l}\tilde{s} = 0; \ \
\tilde{p}\tilde{l} = \tilde{q}\tilde{l} = 0;\nn\\
r\tilde{m} = \tilde{m}r = 0; \ \
r\tilde{n} = \tilde{n}r = 0; \ \
r\tilde{t} = r\tilde{s} = 0; \ \
\tilde{p}r = \tilde{q}r = 0; \nn\\
\tilde{s}\tilde{t} =\tilde{t}\tilde{s} = 0;  \ \
\tilde{s}\tilde{m} = \tilde{s}\tilde{n} = 0; \ \
\tilde{t}\tilde{m} = \tilde{t}\tilde{n} = 0;\nn\\
\tilde{p}\tilde{q} =  \tilde{q}\tilde{p} = 0;\ \
\tilde{m}\tilde{p} = \tilde{n}\tilde{p} = 0; \ \
\tilde{m}\tilde{q} = \tilde{n}\tilde{q} = 0;\nn\\
\tilde{p}^2 =   \tilde{q}^2 = \tilde{s}^2 = \tilde{t}^2= 0.  \end{eqnarray}

The structure of this algebra, denoted ~$\ca_{-+-}\,$, ~is very complicated. There are two free subalgebras:
~$\ca^1_{-+-}$~ with generators ~$\tilde{k},\tilde{l}r$,
~$\ca^2_{-+-}$~ with generators ~$\tilde{m},\tilde{n}$, and four quasi-free subalgebras:
~$\ca^3_{-+-}$~ with generators ~$\tilde{p},\tilde{s}$,
~$\ca^4_{-+-}$~ with generators ~$\tilde{q},\tilde{t}$,
~$\ca^5_{-+-}$~ with generators ~$\tilde{p},\tilde{t}$,
~$\ca^6_{-+-}$~ with generators ~$\tilde{q},\tilde{s}$.
The first four subalgebras have generators as in the ~$\ca_{+-+}$~ case
(but taking into account the last line of \eqref{mpm}).
Only the first two subalgebras are in direct sum, and there are intersections between the last four.
Furthermore, all are related by the following 20 two-letter building blocks:
~$\tilde{k}\tilde{p}\,,\tilde{k}\tilde{q}\,,\tilde{s}\tilde{k},\,\tilde{t}\tilde{k}$,
~$\tilde{l}\tilde{p}\,,\tilde{l}\tilde{q}\,,\tilde{s}\tilde{l},\,\tilde{t}\tilde{l}$,
~$r\tilde{p}\,, r \tilde{q}\,, \tilde{s}r\,,  \tilde{t}r$,
~$\tilde{m}\tilde{s}\,, \tilde{m}\tilde{t}\,,
\tilde{p}\tilde{m}\,, \tilde{q}\tilde{m}$,
~$\tilde{n}\tilde{s}\,, \tilde{n}\tilde{t}\,,
\tilde{p}\tilde{n}\,, \tilde{q}\tilde{n}$, which are the same as in the ~$\ca_{+-+}$~ case,
except those involving ~$r$.

The last difference makes things better. Indeed, there is no overall
ordering, more precisely   we have:
\eqn{ord}\tilde{p},\tilde{q} ~>~ \tilde{m},\tilde{n} ~>~ \tilde{s},\tilde{t} ~>~
\tilde{k},\tilde{l},r ~>~ \tilde{p},\tilde{q}  \ee
i.e., we have some cyclic order.

Thus, the bialgebra ~$\ca_{-+-}$~ may turn out to be the easiest to handle, as the
exotic bialgebra ~$S03$, \cite{ACDM2}.

\bigskip\nt\bu for $a_+ = -a_- = -b_- = c_-$ we have the set of relations
\begin{equation} (-,-,+) = \{ N \cup A_- \cup B_- \cup C_+ \cup
(AB)_+ \cup (AC)_- \cup (BC)_- \} \end{equation}

Explicitly we have:
\begin{eqnarray}
\tilde{k}\tilde{m} =\tilde{m}\tilde{k} = 0;\ \
\tilde{k}\tilde{n} = \tilde{n}\tilde{k} =0; \ \
\tilde{k}\tilde{p} = \tilde{k}\tilde{q} =0; \ \
\tilde{s}\tilde{k} = \tilde{t}\tilde{k} =0;\nn\\
\tilde{l}\tilde{m} = \tilde{m}\tilde{l} = 0;  \ \
\tilde{l}\tilde{n} =  \tilde{n}\tilde{l} = 0; \ \
\tilde{l}\tilde{p} = \tilde{l}\tilde{q} = 0;\ \
\tilde{s}\tilde{l} =\tilde{t}\tilde{l} = 0; \nn\\
r\tilde{m} = \tilde{m}r = 0; \ \
r\tilde{n} = \tilde{n}r = 0; \ \
r\tilde{p} =  r\tilde{q} = 0 ; \ \
\tilde{t}r =  \tilde{s}r = 0; \nn\\
\tilde{s}\tilde{t} =\tilde{t}\tilde{s} = 0;  \ \
\tilde{m}\tilde{s} = \tilde{n}\tilde{s} = 0; \ \
\tilde{m}\tilde{t} = \tilde{n}\tilde{t} = 0;\nn\\
 \tilde{p}\tilde{q} = \tilde{q}\tilde{p} =0; \ \
\tilde{p}\tilde{m} = \tilde{p}\tilde{n} = 0; \ \
\tilde{q}\tilde{m} = \tilde{q}\tilde{n} =  0;\nn\\
\tilde{p}^2 =   \tilde{q}^2 = \tilde{s}^2 = \tilde{t}^2= 0.
  \end{eqnarray}

This algebra, denoted ~$\ca_{--+}\,$, ~is a conjugate of the previous algebra obtained
by the exchange of the pairs of generators ~$(\tilde{p},\tilde{q})$~ and ~$(\tilde{s},\tilde{t})$.

\bigskip\nt\bu for $a_+ = - a_- = - b_- = - c_-$ we have the set of
relations
\begin{equation} (-,-,-) = \{ N \cup A_- \cup B_- \cup C_- \cup
(AB)_+ \cup (AC)_+ \cup (BC)_+ \} \end{equation}

Explicitly we have:
\begin{eqnarray}
\tilde{k}\tilde{m} = \tilde{m}\tilde{k} = 0;\ \
\tilde{k}\tilde{n} = \tilde{n}\tilde{k} =0; \ \
\tilde{k}\tilde{s} = \tilde{s}\tilde{k} = 0;\ \
\tilde{k}\tilde{t} = \tilde{t}\tilde{k} = 0;\nn\\
\tilde{l}\tilde{m} = \tilde{m}\tilde{l} = 0;  \ \
\tilde{l}\tilde{n} =  \tilde{n}\tilde{l} = 0; \ \
\tilde{l}\tilde{s} = \tilde{s}\tilde{l} = 0;\ \
\tilde{l}\tilde{t} = \tilde{t}\tilde{l} = 0;\nn\\
r\tilde{m} = \tilde{m}r = 0; \ \
r\tilde{n} = \tilde{n}r = 0; \ \
r\tilde{p} = \tilde{p}r =  0;\ \
r\tilde{q} =  \tilde{q}r = 0; \nn\\
\tilde{p}\tilde{m} = \tilde{m}\tilde{p} =0;\ \
\tilde{p}\tilde{n} = \tilde{n}\tilde{p} = 0; \ \
\tilde{p}\tilde{s} = \tilde{s}\tilde{p} = 0;\nn\\
\tilde{q}\tilde{m} = \tilde{m}\tilde{q} = 0;\ \
\tilde{q}\tilde{n} = \tilde{n}\tilde{q} = 0;\ \
\tilde{q}\tilde{t} = \tilde{t}\tilde{q} = 0;\nn\\
\tilde{p}^2 = \tilde{t}^2=   \tilde{q}^2 = \tilde{s}^2 = 0 .
\end{eqnarray}

This algebra, denoted ~$\ca_{---}\,$, ~is a conjugate of the algebra ~$\ca_{-++}\,$, obtained
by the exchanges of generators: ~$\tilde{p} \lra\tilde{s}$~ and ~$\tilde{q}\lra\tilde{t}$.

\subsection{Summary}

Thus, taking into account conjugation, we have found ~{\it four}~ different bialgebras:
\eqn{bial}
\ca_{+++} ~\cong~ \ca_{+--} \ , \quad
\ca_{---} ~\cong~ \ca_{-++} \ , \quad
\ca_{+-+} ~\cong~ \ca_{++-} \ , \quad
\ca_{-+-} ~\cong~ \ca_{--+}
\ee
Thus, for future use we shall use shorter notation:
\eqn{bal}
\ca_{++} ~\equiv~ \ca_{+++} \ , \quad
\ca_{--} ~\equiv~ \ca_{---} \ , \quad
\ca_{+-} ~\equiv~ \ca_{+-+} \ , \quad
\ca_{-+} ~\equiv~ \ca_{-+-}
\ee

The first two bialgebras have no ordering. The first one is simpler,
since it is split in two subalgebras with five and four generators.
The third bialgebra has partial ordering in one subalgebra. The last
one, is the most promising since it has partial cyclic ordering.

The next task in this line of research is to find the dual
bialgebras, analogously, as done for the four-element exotic
bialgebras in \cite{ACDM1,ACDM2}. We do this in the next Section for
the most interesting of the above: ~$\ca_{-+} ~\equiv~
\ca_{-+-}\,$.

\section{The dual bialgebra of $A_{-+}$}
\setcounter{equation}{0}

To start with we begin with the coproducts of the elements of the
$T$-matrix. We have:
\begin{equation} \delta { \begin {pmatrix} k & p & l \cr q & r & s
\cr m & t & n \cr \end{pmatrix} } = \begin{pmatrix} k\otimes k
+ p\otimes q + l \otimes m & k\otimes p + p\otimes r + l \otimes t &
k\otimes l + p\otimes s + l\otimes n \cr q\otimes k + r\otimes q +
s\otimes m & q\otimes p + r\otimes r + s \otimes t & q\otimes l +
r\otimes s + s \otimes n \cr m\otimes k + t\otimes q +  n\otimes m &
mo\times p + t\otimes r +  n\otimes t & m\otimes l + t\otimes s +
n\otimes n \cr \end{pmatrix}
\end{equation}
We have also
\begin{equation} \epsilon{\begin {pmatrix} k & p & l \cr q & r & s
\cr m & t & n \cr \end{pmatrix} } = \begin{pmatrix} 1 & 0 & 0 \cr 0 &
1 & 0 \cr 0 & 0 & 1 \cr \end{pmatrix}
\end{equation}
Going to the set of "tilde" generators we see that the co-products
are more complicated. Namely
\begin{eqnarray} \delta(\tilde{k}) = \tilde{k}\otimes\tilde{k} + \tilde{n}\otimes\tilde{n}
+ \tilde{l}\otimes\tilde{l} - \tilde{m}\otimes\tilde{m} +
\tilde{p}\otimes\tilde{q} + \tilde{t}\otimes\tilde{s}, \nn\\
\delta(\tilde{n}) = \tilde{k}\otimes\tilde{n} +
\tilde{n}\otimes\tilde{k} - \tilde{l}\otimes\tilde{m} +
\tilde{m}\otimes\tilde{l} +
\tilde{p}\otimes\tilde{s} + \tilde{t}\otimes\tilde{q}, \nn\\
\delta(\tilde{l}) = \tilde{k}\otimes\tilde{l} +
\tilde{n}\otimes\tilde{m} + \tilde{l}\otimes\tilde{k} -
\tilde{m}\otimes\tilde{n} +
\tilde{p}\otimes\tilde{q} - \tilde{t}\otimes\tilde{s}, \nn\\
\delta(\tilde{m}) = \tilde{k}\otimes\tilde{m} +
\tilde{n}\otimes\tilde{l} - \tilde{l}\otimes\tilde{n} +
\tilde{m}\otimes\tilde{k}
\tilde{p}\otimes\tilde{s} + \tilde{t}\otimes\tilde{q}, \nn\\
\delta(r) = r\otimes r + 2\tilde{q}\otimes\tilde{p} +
2\tilde{s}\otimes\tilde{t}, \nn\\
\delta(\tilde{p}) = \tilde{k}\otimes\tilde{p} +
\tilde{n}\otimes\tilde{t} + \tilde{l}\otimes\tilde{p} -
\tilde{m}\otimes\tilde{t} +
\tilde{p}\otimes r,  \nn\\
\delta(\tilde{t}) = \tilde{k}\otimes\tilde{t} +
\tilde{n}\otimes\tilde{p} - \tilde{l}\otimes\tilde{t} +
\tilde{m}\otimes\tilde{p} +
\tilde{t}\otimes r,  \nn\\
\delta(\tilde{q}) = \tilde{q}\otimes\tilde{k} +
\tilde{q}\otimes\tilde{l} + \tilde{s}\otimes\tilde{n} +
\tilde{s}\otimes\tilde{m} +
r\otimes\tilde{q},  \nn\\
\delta(\tilde{s}) = \tilde{q}\otimes\tilde{n} -
\tilde{q}\otimes\tilde{m} + \tilde{s}\otimes\tilde{k} -
\tilde{s}\otimes\tilde{m} + r\otimes\tilde{s}.
\end{eqnarray}

However, looking at the last four equations we see that they can be
essentially simplified we introduce the new "hat" generators:
\begin{equation} \tilde{k} = \hat{k} + \hat{l}, \ \tilde{l} = \hat{k} -
\hat{l}; \ \ \tilde{m} = \hat{m} - \hat{n}, \ \tilde{n} = \hat{m} +
\hat{n}.  \end{equation} Thus we have:
\begin{eqnarray} \delta(\hat{k}) = 2\hat{k}\otimes\hat{k} -
2\hat{m}\otimes\hat{n} + \tilde{p}\otimes\tilde{q}, \nn\\
\delta(\hat{l}) = 2\hat{l}\otimes\hat{l} -
2\hat{n}\otimes\hat{m} + \tilde{t}\otimes\tilde{s}, \nn\\
 \delta(\hat{m}) = 2\hat{k}\otimes\hat{m} -
2\hat{m}\otimes\hat{l} - \tilde{p}\otimes\tilde{s}, \nn\\
\delta(\hat{n}) = 2\hat{l}\otimes\hat{n} +
2\hat{n}\otimes\hat{k} + \tilde{t}\otimes\tilde{q}, \nn\\
\delta(\hat{p}) = 2\hat{k}\otimes\hat{p} -
2\hat{m}\otimes\hat{t} + \tilde{p}\otimes r, \nn\\
\delta(\hat{q}) = 2\hat{q}\otimes\hat{k} +
2\hat{s}\otimes\hat{n} + r\otimes\tilde{q}, \nn\\
 \delta(\hat{s}) = 2\hat{s}\otimes\hat{l} -
2\hat{q}\otimes\hat{m} + r\otimes\tilde{s}, \nn\\
\delta(\hat{t}) = 2\hat{l}\otimes\hat{t} +
2\hat{n}\otimes\hat{p} + \tilde{t}\otimes r, \nn\\
\delta(r) = r\otimes r + 2\tilde{q}\otimes\tilde{p} +
2\tilde{s}\otimes\tilde{t}.
\end{eqnarray}

With the "hat" generators we have
\begin{eqnarray} \epsilon(\hat{k}) = \epsilon(\hat{l}) = 1/2, \ \
\epsilon(r) = 1, \nn\\ \epsilon(z) = 0, \ {\rm for} \
z=(\hat{m},\hat{n}, \tilde{p},\tilde{q},\tilde{s},\tilde{t}).
\end{eqnarray}

The bialgebra relations are as follows:
\begin{eqnarray}
\hat{k}\hat{m} = \hat{m}\hat{k} = \hat{k}\hat{n} = \hat{n}\hat{k}
=0; \  \hat{k}\tilde{t} = \hat{k}\tilde{s} =
\tilde{p}\hat{k} = \tilde{q}\hat{k} = 0; \nn\\
\hat{l}\hat{m} = \hat{m}\hat{l} = \hat{l}\hat{n} = \hat{n}\hat{l} =
0; \  \hat{l}\tilde{t} = \hat{l}\tilde{s} =
\tilde{p}\hat{l} = \tilde{q}\hat{l} = 0;\nn\\
r\hat{m} = \hat{m}r = r\hat{n} = \hat{n}r = 0; \ \
r\tilde{t} = r\tilde{s} = \tilde{p}r = \tilde{q}r = 0; \nn\\
\tilde{s}\tilde{t} =\tilde{t}\tilde{s} = 0;  \ \ \tilde{s}\hat{m} =
\tilde{s}\hat{n} = \tilde{t}\hat{m} =
\tilde{t}\hat{n} = 0;\nn\\
\tilde{p}\tilde{q} =  \tilde{q}\tilde{p} = 0;\ \ \hat{m}\tilde{p} =
\hat{n}\tilde{p} =
\hat{m}\tilde{q} = \hat{n}\tilde{q} = 0;\nn\\
\tilde{p}^2 =   \tilde{q}^2 = \tilde{s}^2 = \tilde{t}^2= 0.
\end{eqnarray}

The dual elements are defined by the standard procedure of
\cite{Dobrev}, namely:
\begin{equation} \left\langle {Z,f} \right\rangle =
\epsilon\left(\frac {\partial{f}}{\partial{z}} \right), \ \ {\rm
where}\ \ z =
(\hat{k},\hat{l},\hat{m},\hat{n},\tilde{p},\tilde{q},\tilde{s},\tilde{t},r).
\end{equation}

The basis we are working is essentially the following
\begin{eqnarray} \hat{k}^\kappa \hat{l}^\ell r^\tau , \ \ {\rm and \ all \
permutations\ of}\ (\hat{k}\hat{l}r),\nn\\
\hat{k}^\kappa \hat{l}^\ell r^\tau \tilde{p}, \ \ {\rm and \ all \
permutations\
of}\   (\hat{k}\hat{l}r), \nn\\
\hat{k}\hat{l}r)\tilde{q}, \ \ {\rm and \ all \ permutations\ of}\
  (\hat{k}\hat{l}r), \nn\\
\tilde{s}\hat{k}^\kappa \hat{l}^\ell r^\tau , \ \ {\rm and \ all \
permutations\
of}\   (\hat{k}\hat{l}r), \nn\\
\tilde{t}\hat{k}^\kappa \hat{l}^\ell r^\tau  \ \ {\rm and \ all \
permutations\
of}\   (\hat{k}\hat{l}r), \nn\\
\hat{m}, \ \ \hat{n}. \end{eqnarray} Thus the following dual
bialgebra is obtained:
\begin{eqnarray}
\left[\hat{K}^\kappa ,\hat{L}^\ell \right] = 0, \nn\\
\hat{K}^\kappa  \hat{M} = 2^\kappa \hat{M}, \ \
\hat{N}\hat{K}^\kappa  = 2^\kappa  \hat{N},
\ \ \hat{M}\hat{K} = \hat{K}\hat{N} = 0, \nn\\
\hat{M}\hat{L}^\ell  = 2^\ell \hat{M}, \ \ \hat{L}^\ell \hat{N} =
2^\ell  \hat{N},
\ \ \hat{N}\hat{L} = \hat{L}\hat{M} = 0, \nn\\
\hat{M}^2 = \hat{N}^2 = 0, \ \ \hat{M}\hat{N} = - 2 \hat{K}, \ \
\hat{N}\hat{M} = - 2 \hat{L}, \nn\\
\left[\hat{K},\tilde{P}\right] = 2 \tilde{P}, \ \
\left[\hat{L},\tilde{P}\right] = 0,
\ \ \left[\hat{R},\tilde{P}\right] = - \tilde{P}, \nn\\
\left[\hat{K},\tilde{Q}\right] = - 2\tilde{Q}, \ \
\left[\hat{L},\tilde{Q}\right] = 0,
\ \ \left[\hat{R},\tilde{Q}\right] =  \tilde{Q}, \nn\\
\left[\hat{K},\tilde{S}\right] = 0, \ \
\left[\hat{L},\tilde{S}\right] = - 2\tilde{S},
\ \ \left[\hat{R},\tilde{S}\right] =  \tilde{S}, \nn\\
\left[\hat{K},\tilde{T}\right] = 0, \ \
\left[\hat{L},\tilde{T}\right] = 2\tilde{T},
\ \ \left[\hat{R},\tilde{T}\right] = - \tilde{T}, \nn\\
\hat{M}\tilde{T} = - 2\tilde{P}, \ \ \tilde{T}\hat{M} = 0, \ \
\tilde{Q}\hat{M} = - 2\tilde{S}, \ \ \hat{M}\tilde{Q} = 0, \nn\\
\hat{N}\tilde{P} = 2\tilde{T}, \ \ \tilde{P}\hat{N} = 0, \ \
\tilde{S}\hat{N} = 2\tilde{Q}, \ \ \hat{N}\tilde{S} = 0, \nn\\
\hat{M}\tilde{P} = \tilde{P}\hat{M} = \hat{M}\tilde{S} =
\tilde{S}\hat{M} = 0, \nn\\
\hat{N}\tilde{Q} = \tilde{Q}\hat{N} = \hat{N}\tilde{T} =
\tilde{T}\hat{N} = 0, \nn\\
\left[\tilde{S},\tilde{P}\right] = \hat{M}, \ \
\left[\tilde{Q},\tilde{T}\right] = \hat{N}, \ \
\left[\tilde{Q},\tilde{S}\right] = \left[\tilde{P},\tilde{T}\right]
=
0, \nn\\
\tilde{P}\tilde{Q} = \tilde{T}\tilde{S}, \ \ \tilde{Q}\tilde{P} =
\tilde{S}\tilde{T}, \nn\\
\tilde{P}^2 = \tilde{Q}^2 = \tilde{S}^2 = \tilde{T}^2 = 0, \ \
\tilde{P}\tilde{T} = \tilde{T}\tilde{P} = \tilde{Q}\tilde{S} =
\tilde{S}\tilde{Q} = 0.
\end{eqnarray}
Finally we write down the co-products of the dual bialgebra:
\begin{eqnarray} \delta(\hat{K}) = \hat{K}\otimes 1_U + 1_U\otimes
\hat{K}, \nn\\ \delta(\hat{L}) = \hat{L}\otimes 1_U + 1_U\otimes
\hat{L}, \nn\\ \delta(\hat{M}) = \hat{M}\otimes 1_U + 1_U\otimes
\hat{M}, \nn\\ \delta(\hat{N}) = \hat{N}\otimes 1_U + 1_U\otimes
\hat{N}, \nn\\ \delta(\tilde{P}) = \tilde{P}\otimes 1_U + 1_U\otimes
\tilde{P}, \nn\\ \delta(\tilde{Q}) = \tilde{Q}\otimes 1_U + 1_U\otimes
\tilde{Q}, \nn\\ \delta(\tilde{S}) = \tilde{S}\otimes 1_U +
1_U\otimes\tilde{S}, \nn\\ \delta(\tilde{T}) = \tilde{T}\otimes 1_U +
1_U\otimes\tilde{T}, \nn\\ \delta(R) = R\otimes 1_U + 1_U\otimes R.
\end{eqnarray}.
%%%%%%

\section{Summary and outlook}

In the present paper we have found a multitude of exotic bialgebras
and the dual  of  one of them. More duals should be constructed. More
importantly, we plan to continue the programme fulfilled
successfully for the exotic bialgebra S03, cf. summary in
\cite{ACDM4}. In particular, it is important to find the FRT duals
\cite{FRT} which are different from the standard duals for the
exotic bialgebras. Further we should find the Baxterisation of the
dual algebras.   Their finite-dimensional representations  should be
considered. Diagonalisations of the braid matrices would be used to
handle  the representations of the corresponding $L$-algebras (in
the FRT formalism) and to formulate the fusion of finite-dimensional
representations. Possible applications should be considered, in
particular,   exotic  vertex models and   integrable spin-chain
models.

\section*{Acknowledgments}

One of the authors (VKD) would like to thank for hospitality the
Abdus Salam International Center for Theoretical Physics, where part
of the work was done. VKD and SGM were supported in part by
Bulgarian NSF grant {\it DO 02-257}.

%\section*{References}

\end{document}